\documentclass[preprint,11pt]{elsarticle}

\usepackage[utf8]{inputenc}
\usepackage[T1]{fontenc}

\usepackage{amsmath}
\usepackage{amssymb}
\usepackage{amsthm}     
\usepackage{dsfont}     
\usepackage{gensymb}    
\usepackage{textcomp}   

\usepackage{graphicx}
\usepackage{subcaption} 
\usepackage{xcolor}

\usepackage[colorlinks=true, linkcolor=blue, citecolor=blue, urlcolor=blue]{hyperref}

\usepackage{pgf}
\usepackage{tikz}
\usetikzlibrary{decorations.text,positioning,fit,patterns,shapes, arrows,calc,decorations.pathreplacing,trees}
\usepackage{pgfplots}  
\usepackage{tikz-network}
\usepackage{circuitikz}

\newtheorem{remark}{Remark}

\usepackage{mathtools}
\mathtoolsset{showonlyrefs}

\begin{document}

\begin{frontmatter}

\title{Physics-Informed Deep Learning for Industrial Processes: Time-Discrete VPINNs for heat conduction}


\author[1]{Manuela Bastidas Olivares\corref{cor1}}
\ead{mbastidaso@unal.edu.co}

\author[1]{Josué David Acosta Castrillón}

\author[1]{Diego A. Muñoz}


\affiliation[1]{
    organization={Universidad Nacional de Colombia, Sede Medellín},
    addressline={Facultad de Ciencias}, 
    city={Medellín},
    country={Colombia}
}

\begin{abstract}
Neural networks offer powerful tools to solve partial differential equations (PDEs). We present a Variational Physics-Informed Neural Network (VPINN) designed for parabolic problems. Our approach combines a classical time discretization with a composed loss function, which minimizes the residual's dual norm at every time step. We validate the framework by modeling the freezing of coffee extracts in an industrial cylinder. The simulation accounts for temperature-dependent properties and experimental data. It successfully captures the thermal dynamics of the process. 
\end{abstract}


\end{frontmatter}


\section{Introduction}
\label{sec:intro}

Parabolic partial differential equations (PDEs) are fundamental to modeling diffusion processes in science and engineering, with the heat equation serving as the primary prototype. Then, developing efficient numerical methods for these types of equations remains important and even more now when Physics-Informed Neural Networks (PINNs) have emerged as a powerful, mesh-free alternative to classical solvers. By embedding the differential operator directly into the loss function, the classical PINNs solve PDEs without complex grid generation and computational advantages, as mentioned in \cite{cuomo2022scientific, lu2021deepxde, raissi2019physics}. Recent literature extensively documents their application to parabolic problems \cite{ beck2019machine, cai2021physics, meng2020ppinn, nusken2023interpolating, shin2020convergence}. However, standard PINNs typically target the strong form of the equation, which demands high regularity of the solution, but this constraint becomes a limitation when solutions lack smoothness or exhibit sharp gradients, as in usual in real-life applications \cite{krishnapriyan2021characterizing,wang2021understanding,wang2022and}.

To overcome these theoretical restrictions, we adopt the weak variational formulation, as in the usual finite element implementations, transferring derivatives to (regular enough) test functions, allowing for solutions with lower regularity. We propose a robust Variational PINN (RVPINN) framework that integrates variational principles with a time discretization. Instead of using the mean-squared error \cite{berrone2022solving, kharazmi2019variational,kharazmi2021hp}, our method minimizes the residual's dual norm at each time step. This norm acts as a rigorous error estimator, as shown in \cite{rojas2024robust,taylor2023deep} and used in \cite{taylor2025deep, uriarte2025optimizing}.

We demonstrate the practical value of this framework in the context of industrial freezing. In particular, we model the heat transfer processes in coffee extracts, which are characterized by temperature-dependent thermophysical properties. Here, we use a  reference dataset to characterize the material properties: density, thermal conductivity, and heat capacity—across the phase change temperature range and the experimental records provide the necessary boundary conditions. To ensure physical fidelity, our main numerical experiment is validated with experimental data obtained from an industrial coffee extract freezing process. The specific experimental setup, measurement protocols, and raw data acquisition appear in \cite{restrepo2025container}. Based on this experimental framework, we confirm that the method accurately captures complex thermal behaviors in real-world scenarios.

This paper is organized as follows. Section \ref{sec:formulation} details the mathematical framework, deriving the time-discrete variational formulation for parabolic problems and a classical well-posedness analysis. In Section \ref{ConstrucctionLoss}, we introduce the concept of a neural network as an approximation to the solution of a PDE. There, we construct a loss function that accurate approximate the residual's dual norm as an estimator of the error.  We validate the proposed method in Section \ref{sec:experiments} through numerical experiments on a benchmark heat equation and the industrial coffee extract application. Finally, Section \ref{sec:conclusions} summarizes our findings and discusses future directions.

\section{Mathematical Framework}
\label{sec:formulation}

We consider a time-dependent diffusion problem defined on an open bounded domain $\Omega \subset \mathbb{R}^d$ with a Lipschitz boundary $\partial \Omega$. The physical process is governed by the conservation of energy, resulting in a non-linear parabolic partial differential equation with solution-dependent coefficients.

Let $u(x,t)$ denote the scalar field of interest (e.g., temperature) for $x \in \Omega$ and $t \in (0, t_{\text{end}}]$. The governing equation is given by:
\begin{equation}\label{eq:strong_form}
    \frac{\partial}{\partial t} \Big( \mathcal{C}(u)\, u \Big) \;-\; \nabla \cdot \Big( \mathcal{K}(u)\, \nabla u \Big) \;=\; f(x,t),
\end{equation}
subject to the boundary and initial conditions:
\begin{align*}
    u(x,t) &= u_b(t), \quad && \text{on } \partial \Omega \times (0, t_{\text{end}}], \\
    u(x,0) &= u_I(x), \quad && \text{in } \Omega.
\end{align*}

In this energy balance model, the internal energy term is represented by the product $\mathcal{C}(u)u$, where the coefficient $\mathcal{C}(\cdot)$  corresponds to the volumetric heat capacity, defined as the product of the temperature-dependent specific heat and the density of the material. The diffusive contribution to the balance is given by $\mathcal{K} \nabla u$, representing heat transfer by conduction. This specific mechanism is modeled by Fourier's Law, where $\mathcal{K}(\cdot)$ serves as the thermal conductivity or diffusion coefficient. We assume that all the parameters are strictly positive and bounded functions of the state variable $u$. The term $f(x,t)$ represents a volumetric source term. Moreover, $u_I$ is the initial condition, and $u_b$ is the boundary condition. 

\begin{remark}
    For the sake of simplicity, we assume homogeneous Dirichlet boundary conditions and consider the linear case, where the model parameters depend solely on the spatial coordinates rather than on $u$. These assumptions allow for a clearer exposition of the mathematical framework. Subsequently, in Section~\ref{sec:experiments}, we extend the method to the general nonlinear case, demonstrating robustness in more complex scenarios.
\end{remark}

\subsection{Time Discretization}

We discretize the time interval $(0, t_{\text{end}}]$ into $N$ steps of size $\Delta t$, denoting the time at step $n$ as $t^n = n\Delta t$ and the solution as $u^n(x) \approx u(x, t^n)$.

Approximating the time derivative with a backward Euler discretization, we obtain the semi-discrete equation:
\begin{equation*}
    \mathcal{C} \cdot \frac{ u^n \;-\; u^{n-1}}{\Delta t} \;-\; \nabla \cdot \Big( \mathcal{K} \, \nabla u^n \Big) \;=\; f(x,t^n).
\end{equation*}

Rearranging terms, we arrive at a linear elliptic problem to be solved sequentially for each $n$:
\begin{equation}\label{eq:semi_discrete_linear}
    \mathcal{C}\, u^n \;-\; \Delta t\, \nabla \cdot \Big( \mathcal{K}\, \nabla u^n \Big) \;=\; \underset{:= f^{n}(x)}{\underbrace{\mathcal{C} \, u^{n-1} + \Delta t\,f(x,t^n)}}.
\end{equation}

\subsection{Variational Formulation}

To derive the weak formulation, we multiply the semi-discrete equation \eqref{eq:semi_discrete_linear} by a test function $v \in H_0^1(\Omega)$ and integrating over the domain $\Omega$ (applying Green's first identity), we have:
\begin{equation*}
    \begin{split}
        - \int_{\Omega} \nabla \cdot \Big( \mathcal{K}\, \nabla u^n \Big)\, v \, dx \;&=\; \int_{\Omega} \mathcal{K}\, \nabla u^n \cdot \nabla v \, dx \\  & \qquad -\; \int_{\partial \Omega} \mathcal{K} (\nabla u^n \cdot \mathbf{n})\, v \, ds,
    \end{split}
\end{equation*}
where $\mathbf{n}$ is the outward unit normal vector to $\partial \Omega$. Since $v \in H_0^1(\Omega)$, the trace of $v$ on the boundary vanishes and the equation reduces to:
\[
    \int_{\Omega} \mathcal{C}\, u^n\, v \, dx \;+\; \Delta t \int_{\Omega} \mathcal{K}\, \nabla u^n \cdot \nabla v \, dx \;=\; \int_{\Omega} f^{n}\, v \, dx.
\]

Thus, the variational problem at time step $n$ is to find $u^n \in H_0^1(\Omega)$ such that:
\begin{equation}\label{eq:variational_general}
    a(u^n, v) \;=\; l^n(v), \quad \forall v \in H_0^1(\Omega),
\end{equation}
where the parameter-dependent bilinear form $a(\cdot, \cdot)$ and linear functional $l(\cdot)$ are defined as:
\begin{align}
    a(w, v) &:=\; \int_{\Omega} \mathcal{C}\, w\, v \, dx \;+\; \Delta t \int_{\Omega} \mathcal{K}\, \nabla w \cdot \nabla v \, dx, \label{eq:bilinear_def} \\[6pt]
    l^n(v) &:=\; \int_{\Omega} f^{n}  \, v \, dx. \label{eq:linear_def}
\end{align}


\subsection{Well-Posedness Analysis}

The solvability of the sequence of variational problems is guaranteed by the Lax-Milgram theorem. To apply this theorem, we must verify that the bilinear form $a(\cdot, \cdot)$ is both continuous and coercive on $H_0^1(\Omega) \times H_0^1(\Omega)$.

We equip the space $H_0^1(\Omega)$ with the gradient norm:
\[
    \|v\|_{H_0^1} := \|\nabla v\|_{L^2(\Omega)}.
\]
Due to the Poincaré inequality, there exists a constant $C_P > 0$ (depending only on $\Omega$) such that $\|v\|_{L^2} \le C_P \|v\|_{H_0^1}$ for all $v \in H_0^1(\Omega)$. We assume physically consistent material properties, meaning there exist positive constants such that:
\[
    0 < c_{\min} \le \mathcal{C}(\cdot) \le c_{\max}, \quad \text{and} \quad 0 < k_{\min} \le \mathcal{K}(\cdot) \le k_{\max}.
\]

\textbf{Continuity:} Using the Cauchy-Schwarz inequality and the upper bounds of the coefficients, we have:
\begin{align*}
    |a(w, v)| &\le \int_{\Omega} |\mathcal{C}\, w\, v| \, dx + \Delta t \int_{\Omega} |\mathcal{K}\, \nabla w \cdot \nabla v| \, dx \\
    &\le c_{\max} \|w\|_{L^2} \|v\|_{L^2} + \Delta t\, k_{\max} \|\nabla w\|_{L^2} \|\nabla v\|_{L^2}.
\end{align*}
Applying the Poincaré inequality to the $L^2$ terms, we obtain:
\[
    |a(w, v)| \le c_{\max} C_P^2 \|w\|_{H_0^1} \|v\|_{H_0^1} + \Delta t\, k_{\max} \|w\|_{H_0^1} \|v\|_{H_0^1}.
\]
Thus, the bilinear form is bounded with continuity constant $M = c_{\max} C_P^2 + \Delta t\, k_{\max}$.

\textbf{Coercivity:} Using the lower bounds of the coefficients, we examine the quadratic form:
\begin{align*}
    a(v, v) &= \int_{\Omega} \mathcal{C}\, v^2 \, dx + \Delta t \int_{\Omega} \mathcal{K}\, |\nabla v|^2 \, dx \\
    &\ge c_{\min} \|v\|_{L^2}^2 + \Delta t\, k_{\min} \|\nabla v\|_{L^2}^2.
\end{align*}
Since the first term $c_{\min} \|v\|_{L^2}^2$ is non-negative, we can bound the expression from below:
\[
    a(v, v) \ge \Delta t\, k_{\min} \|\nabla v\|_{L^2}^2 = \Delta t\, k_{\min} \|v\|_{H_0^1}^2.
\]
Defining $\gamma := \Delta t\, k_{\min}$, we have $a(v, v) \ge \gamma \|v\|_{H_0^1}^2$. Since $\Delta t > 0$ and $k_{\min} > 0$, $\gamma$ is strictly positive, so the form is coercive. By the Lax-Milgram theorem, a unique solution $u^n$ exists for every time step $n$.

Consequently, at each time step, the resulting problem is well-posed, providing a consistent mathematical framework to analyze the approximation errors of the solution.

\subsection{Stability Estimates and Residual-Based Error Bounds}

In the context of the weak formulation of this type of PDEs, the stability of the solution is governed by the \textit{inf-sup} condition (or the Babuška-Brezzi condition). Since our trial and test spaces coincide in the single Hilbert space $V = H_0^1(\Omega)$, the coercivity property established in the previous section is sufficient to guarantee this condition.

Let $\gamma$ be the coercivity constant derived in the previous section. By choosing the test function $v = w$ in the supremum, we obtain:
\[
  \sup_{v \neq 0} \frac{|a(w,v)|}{\|v\|_{V}} 
  \;\ge\; 
  \frac{a(w,w)}{\|w\|_{V}} 
  \;\ge\; 
  \gamma \|w\|_{V}.
\]
Consequently, the bilinear form satisfies the inf-sup condition:
\begin{equation}\label{eq:inf_sup}
  \inf_{w \neq 0} \sup_{v \neq 0}
  \frac{|a(w,v)|}{\|w\|_{V}\,\|v\|_{V}}
  \;\ge\;
  \gamma \;>\; 0.
\end{equation}

The inequality \eqref{eq:inf_sup}, combined with the continuity of the bilinear form (with constant $M$), provides the theoretical foundation for our residual minimization (neural network training) strategy.

Let $u^* \in V$ be the exact solution to the variational problem \eqref{eq:variational_general} at a specific time step $n$, given the previous state $u^{n-1}$. By definition, it satisfies:
\[
    a(u^*, v) = l^n(v), \quad \forall v \in V.
\]
For any approximate solution $u^n \in V$ (e.g., the output of the neural network), we define the residual functional $\mathcal{R}(u^n) \in V^*$ (the dual space $H^{-1}(\Omega)$) as:
\begin{equation}\label{eq:residual}
    \langle \mathcal{R}(u^n), v \rangle_{V^* \times V} := a(u^n, v) - l^n(v).
\end{equation}
    
Due to the linearity of $a(\cdot, \cdot)$ with respect to its first argument, we can rewrite the residual in terms of the error $e = u^n - u^*$:
\[
    \langle \mathcal{R}(u^n), v \rangle_{V^* \times V} \;=\; a(u^n, v) - a(u^*, v) \;=\; a(u^n - u^*, v).
\]

Using the continuity and inf-sup conditions, we can bound the error $\|u^n - u^*\|_V$ in terms of the dual norm of the residual:
\begin{equation}\label{eq:error_estimator}
    \frac{1}{M} \|\mathcal{R}(u^n)\|_{V^*} \;\le\; \|u^n - u^*\|_{V} \;\le\; \frac{1}{\gamma} \|\mathcal{R}(u^n)\|_{V^*}.
\end{equation}

Equation \eqref{eq:error_estimator} establishes that the energy norm of the error is equivalent to the dual norm of the residual. Therefore, minimizing $\|\mathcal{R}(u)\|_{V^*}$ is mathematically equivalent to decreasing the true error of the solution at each time step. 

Now, the primary objective of this work is to design a strategy for the accurate computation of this dual norm at each time step.

\section{Neural Network Discretization}\label{ConstrucctionLoss}
\label{sec:DFR0}

Consistent with our time-stepping formulation, we partition the interval $(0, t_{\text{end}}]$ into $N$ uniform steps and approximate the sequence of spatial solutions $\{u^n\}_{n=1}^{N}$ using a single fully-connected feed-forward neural network. Unlike standard approaches that treat time as an input coordinate, we adopt a discretization strategy where the network maps the spatial coordinate $x \in \Omega$ to a vector of solutions at discrete time steps.

Let $u_\theta: \Omega \to \mathbb{R}^{N}$ denote the neural network parameterized by the set of weights and biases $\theta$. The architecture consists of $L$ layers. The output of the $j$-th layer, denoted by $z_j$, is defined recursively as:
\begin{equation}\label{eq:layerspar}
    z_j = \sigma_j(W_j z_{j-1} + b_j), \quad \text{for } j = 1, \dots, L,
\end{equation}
where $z_0 = x$ is the spatial input. The functions $\sigma_j$ represent non-linear activation functions (e.g., $\tanh$) for the hidden layers, while the output activation $\sigma_L$ is the identity. The final output vector corresponds to the approximation at each time step $ \begin{bmatrix} \hat{u}^1(x), \dots, \hat{u}^N(x) \end{bmatrix}^\top$.

To strictly enforce homogeneous Dirichlet boundary conditions, we employ a non-trainable cutoff function $\chi: \overline{\Omega} \to \mathbb{R}$. This function is constructed to satisfy $\chi(x) = 0$ for all $x \in \partial \Omega$ and $\chi(x) > 0$ in the interior. The final trial solution for the $n$-th time step is given by:
\begin{equation}
    u^n_\theta(x) := \chi(x) \, \hat{u}^n(x),
\end{equation}
and with this, it satisfies the homogeneous boundary conditions.

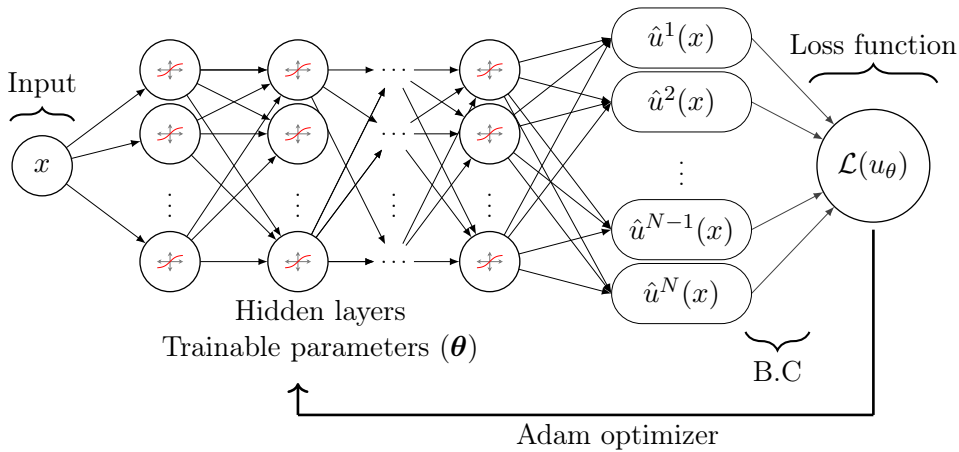
\begin{figure}[htpb!]
    \centering
    \begin{tikzpicture}[scale=0.85]
	
    \draw[decorate, decoration={brace, amplitude=5pt}, thick] (-0.5, 2.2) -- (0.5, 2.2) node[midway, above=5pt] {Input};
    
    \draw (4.35,-1.1) node[align=center] {Hidden layers \\ Trainable parameters ($\boldsymbol{\theta}$)};
	
	\draw[decorate, decoration={brace, amplitude=7pt}, thick] (12, 2.8) -- (14, 2.8) node[midway, above=7pt] {Loss function};
	
	\draw[color=black,  line width=1pt] (13,0.5) -- (13,-2.4) ;
	\draw[color=black,  line width=1pt] (13,-2.4) -- (4,-2.4);
	\draw[->, color=black,  line width=1pt] (4,-2.4) -- (4,-1.9);
	\draw (9,-2.75) node {Adam optimizer};

    \draw[decorate, decoration={brace, amplitude=7pt,mirror}, thick] (11, -1.1) -- (12, -1.1) node[midway, below=7pt] {B.C};
	\Vertex[x=0,y=1.5,label=$x$,color=white,size=0.8,fontsize=\normalsize,style={line width=0.5pt}]{X}
	\Vertex[x=2,y=3,color=white,size=0.8, style={line width=0.5pt}]{A11}
	\Vertex[x=4,y=3,color=white,size=0.8, style={line width=0.5pt}]{A21}
	\Vertex[x=7,y=3,color=white,size=0.8, style={line width=0.5pt}]{A31}
	\Vertex[x=2,y=2,color=white,size=0.8, style={line width=0.5pt}]{A12}
	\Vertex[x=4,y=2,color=white,size=0.8, style={line width=0.5pt}]{A22}
	\Vertex[x=7,y=2,color=white,size=0.8, style={line width=0.5pt}]{A32}
	\Vertex[x=2,y=0,color=white,size=0.8, style={line width=0.5pt}]{A13}
	\Vertex[x=4,y=0,color=white,size=0.8, style={line width=0.5pt}]{A23}
	\Vertex[x=7,y=0,color=white,size=0.8, style={line width=0.5pt}]{A33}
	
	\begin{scope}[shift={(A11.center)}, scale=0.1]
		\draw[color=gray,<->,>={Stealth[length=0.5mm,width=0.5mm]}](-2,0) -- (2,0) ;
		\draw[color=gray,<->,>={Stealth[length=0.5mm,width=0.5mm]}](0,-1.5) -- (0,1.5) ;
		\draw[domain=-2:2,smooth,variable=\x,red,line width=0.3pt] plot ({\x},{tanh(\x)});
	\end{scope}
	\begin{scope}[shift={(A21.center)}, scale=0.1]
		\draw[color=gray,<->,>={Stealth[length=0.5mm,width=0.5mm]}](-2,0) -- (2,0) ;
		\draw[color=gray,<->,>={Stealth[length=0.5mm,width=0.5mm]}](0,-1.5) -- (0,1.5) ;
		\draw[domain=-2:2,smooth,variable=\x,red,line width=0.3pt] plot ({\x},{tanh(\x)});
	\end{scope}
	\begin{scope}[shift={(A31.center)}, scale=0.1]
		\draw[color=gray,<->,>={Stealth[length=0.5mm,width=0.5mm]}](-2,0) -- (2,0) ;
		\draw[color=gray,<->,>={Stealth[length=0.5mm,width=0.5mm]}](0,-1.5) -- (0,1.5) ;
		\draw[domain=-2:2,smooth,variable=\x,red,line width=0.3pt] plot ({\x},{tanh(\x)});
	\end{scope}
	\begin{scope}[shift={(A12.center)}, scale=0.1]
		\draw[color=gray,<->,>={Stealth[length=0.5mm,width=0.5mm]}](-2,0) -- (2,0) ;
		\draw[color=gray,<->,>={Stealth[length=0.5mm,width=0.5mm]}](0,-1.5) -- (0,1.5) ;
		\draw[domain=-2:2,smooth,variable=\x,red,line width=0.3pt] plot ({\x},{tanh(\x)});
	\end{scope}
	\begin{scope}[shift={(A13.center)}, scale=0.1]
		\draw[color=gray,<->,>={Stealth[length=0.5mm,width=0.5mm]}](-2,0) -- (2,0) ;
		\draw[color=gray,<->,>={Stealth[length=0.5mm,width=0.5mm]}](0,-1.5) -- (0,1.5) ;
		\draw[domain=-2:2,smooth,variable=\x,red,line width=0.3pt] plot ({\x},{tanh(\x)});
	\end{scope}
	\begin{scope}[shift={(A22.center)}, scale=0.1]
		\draw[color=gray,<->,>={Stealth[length=0.5mm,width=0.5mm]}](-2,0) -- (2,0) ;
		\draw[color=gray,<->,>={Stealth[length=0.5mm,width=0.5mm]}](0,-1.5) -- (0,1.5) ;
		\draw[domain=-2:2,smooth,variable=\x,red,line width=0.3pt] plot ({\x},{tanh(\x)});
	\end{scope}
	\begin{scope}[shift={(A23.center)}, scale=0.1]
		\draw[color=gray,<->,>={Stealth[length=0.5mm,width=0.5mm]}](-2,0) -- (2,0) ;
		\draw[color=gray,<->,>={Stealth[length=0.5mm,width=0.5mm]}](0,-1.5) -- (0,1.5) ;
		\draw[domain=-2:2,smooth,variable=\x,red,line width=0.3pt] plot ({\x},{tanh(\x)});
	\end{scope}
	\begin{scope}[shift={(A32.center)}, scale=0.1]
		\draw[color=gray,<->,>={Stealth[length=0.5mm,width=0.5mm]}](-2,0) -- (2,0) ;
		\draw[color=gray,<->,>={Stealth[length=0.5mm,width=0.5mm]}](0,-1.5) -- (0,1.5) ;
		\draw[domain=-2:2,smooth,variable=\x,red,line width=0.3pt] plot ({\x},{tanh(\x)});
	\end{scope}
	\begin{scope}[shift={(A33.center)}, scale=0.1]
		\draw[color=gray,<->,>={Stealth[length=0.5mm,width=0.5mm]}](-2,0) -- (2,0) ;
		\draw[color=gray,<->,>={Stealth[length=0.5mm,width=0.5mm]}](0,-1.5) -- (0,1.5) ;
		\draw[domain=-2:2,smooth,variable=\x,red,line width=0.3pt] plot ({\x},{tanh(\x)});
	\end{scope}

	\Vertex[x=10,y=3.5,label={$\hat{u}^1(x)$},color=white,shape=rectangle,fontsize=\normalsize,style={draw=black, rounded rectangle, line width=0.3pt, inner sep=2pt,  minimum width=2cm, minimum height=0.8cm, align=center}]{tildeu1}
	
	\Vertex[x=10,y=2.5,label={$\hat{u}^2(x)$},color=white,shape=rectangle,fontsize=\normalsize,style={draw=black, rounded rectangle, line width=0.3pt, inner sep=2pt,  minimum width=2cm, minimum height=0.8cm, align=center}]{tildeu2}
		
	\Vertex[x=10,y=0.5,label={$\hat{u}^{N-1}(x)$},color=white,shape=rectangle,fontsize=\normalsize,style={draw=black, rounded rectangle, line width=0.3pt, inner sep=2pt,  minimum width=2cm, minimum height=0.8cm, align=center}]{tildeuN}
	
	\Edge[color=white,label={$\color{black}\vdots$}](tildeu2)(tildeuN)
	
	\Vertex[x=10,y=-0.5,label={$\hat{u}^N(x)$},color=white,shape=rectangle,fontsize=\normalsize,style={draw=black, rounded rectangle, line width=0.3pt, inner sep=2pt,  minimum width=2cm, minimum height=0.8cm, align=center}]{tildeuNN}
		
	\Vertex[x=13,y=1.5,label=$\mathcal{L}(u_{\theta})$,color=white,size=1.5,fontsize=\normalsize,style={line width=0.5pt}]{lu}
	
	\Vertex[x=5.5,y=3,style={color=white},label=$\color{black}\hdots$,size=0.5]{Ah1}
	\Vertex[x=5.5,y=2,style={color=white},label=$\color{black}\hdots$,size=0.5]{Ah2}
	\Vertex[x=5.5,y=0,style={color=white},label=$\color{black}\hdots$,size=0.5]{Ah3}

	\Edge[color=black, Direct,lw=0.3pt](X)(A11)
	\Edge[color=black, Direct,lw=0.3pt](X)(A12)
	\Edge[color=black, Direct,lw=0.3pt](X)(A13)
	
	\Edge[color=white,label={$\color{black}\vdots$}](A12)(A13)
	\Edge[color=white,label={$\color{black}\vdots$}](A22)(A23)
	\Edge[color=white,label={$\color{black}\vdots$}](Ah2)(Ah3)
	\Edge[color=white,label={$\color{black}\vdots$}](A32)(A33)
	
	\Edge[color=black, Direct,lw=0.25pt](A11)(A21)
	\Edge[color=black, Direct,lw=0.25pt](A11)(A22)
	\Edge[color=black, Direct,lw=0.25pt](A11)(A23)
	\Edge[color=black, Direct,lw=0.25pt](A11)(A21)
	\Edge[color=black, Direct,lw=0.25pt](A12)(A21)
	\Edge[color=black, Direct,lw=0.25pt](A12)(A22)
	\Edge[color=black, Direct,lw=0.25pt](A12)(A23)
	\Edge[color=black, Direct,lw=0.25pt](A13)(A21)
	\Edge[color=black, Direct,lw=0.25pt](A13)(A22)
	\Edge[color=black, Direct,lw=0.25pt](A13)(A23)
	
	\Edge[color=black, Direct,lw=0.25pt](A23)(Ah1)
	\Edge[color=black, Direct,lw=0.25pt](A23)(Ah2)
	\Edge[color=black, Direct,lw=0.25pt](A23)(Ah3)
	\Edge[color=black, Direct,lw=0.25pt](A21)(Ah1)
	\Edge[color=black, Direct,lw=0.25pt](A21)(Ah2)
	\Edge[color=black, Direct,lw=0.25pt](A21)(Ah3)
	\Edge[color=black, Direct,lw=0.25pt](A23)(Ah1)
	\Edge[color=black, Direct,lw=0.25pt](A23)(Ah2)
	\Edge[color=black, Direct,lw=0.3pt](A23)(Ah3)
	
	\Edge[color=black, Direct,lw=0.3pt](Ah1)(A31)
	\Edge[color=black, Direct,lw=0.3pt](Ah1)(A32)
	\Edge[color=black, Direct,lw=0.3pt](Ah1)(A33)
	\Edge[color=black, Direct,lw=0.3pt](Ah2)(A31)
	\Edge[color=black, Direct,lw=0.3pt](Ah2)(A32)
	\Edge[color=black, Direct,lw=0.3pt](Ah2)(A33)
	\Edge[color=black, Direct,lw=0.3pt](Ah3)(A31)
	\Edge[color=black, Direct,lw=0.3pt](Ah3)(A32)
	\Edge[color=black, Direct,lw=0.3pt](Ah3)(A33)
	
	\Edge[color=black, Direct,lw=0.3pt](A31)(tildeu1.west)
	\Edge[color=black, Direct,lw=0.3pt](A32)(tildeu1.west)
	\Edge[color=black, Direct,lw=0.3pt](A33)(tildeu1.west)
	
	\Edge[color=black, Direct,lw=0.3pt](A31)(tildeu2.west)
	\Edge[color=black, Direct,lw=0.3pt](A32)(tildeu2.west)
	\Edge[color=black, Direct,lw=0.3pt](A33)(tildeu2.west)
	
	\Edge[color=black, Direct,lw=0.3pt](A31)(tildeuN.west)
	\Edge[color=black, Direct,lw=0.3pt](A32)(tildeuN.west)
	\Edge[color=black, Direct,lw=0.3pt](A33)(tildeuN.west)
	
	\Edge[color=black, Direct,lw=0.3pt](A31)(tildeuNN.west)
	\Edge[color=black, Direct,lw=0.3pt](A32)(tildeuNN.west)
	\Edge[color=black, Direct,lw=0.3pt](A33)(tildeuNN.west)
	
	\Edge[Direct,lw=0.3pt](tildeu1.east)(lu)
	\Edge[Direct,lw=0.3pt](tildeu2.east)(lu)
	\Edge[Direct,lw=0.3pt](tildeuN.east)(lu)
	\Edge[Direct,lw=0.3pt](tildeuNN.east)(lu)
	
 \end{tikzpicture}
    \vspace{-0.5cm}
    \caption{Neural network architecture}
    \label{fig:placeholder}
\end{figure}

We determine the network parameters ${\theta}$ by minimizing the total loss $\mathcal{L}(u_{\theta})$ using the {Adam} optimizer with a fixed learning rate. Here, the gradients are computed via automatic differentiation to iteratively (via backpropagation) update the weights and biases of the network. This architecture and optimization process is sketched in Figure~\ref{fig:placeholder}

\subsection{The loss function calculation}
Our objective is to construct a loss function $\mathcal{L}(u_{\theta})$ whose minimization is mathematically equivalent to minimizing the true approximation error.

From the analysis in Section~\ref{sec:formulation}, at each time step the error is $\|u^n - u^*\|_{H_0^1}$, with $u^*$ being the exact solution at the time $t= n\Delta t$. From \eqref{eq:error_estimator}, this error is bounded by the dual norm of the residual $\|\mathcal{R}(u^n)\|_{H^{-1}}$.
By the Riesz Representation Theorem, there exists a unique representative $\psi^n \in H_0^1(\Omega)$ such that:
\[
    \langle \mathcal{R}(u^n), v \rangle = (\psi^n, v)_{H_0^1}, \quad \text{and} \quad \|\mathcal{R}(u^n)\|_{H^{-1}} = \|\psi^n\|_{H_0^1}.
\]
We approximate this norm by projecting the representative onto an orthonormal basis $\{\varphi_k\}_{k=1}^{N_{\text{test}}}$ of the test space $H_0^1(\Omega)$. Using Parseval's identity, the squared norm is given by the sum of the squared coefficients:
\begin{equation}\label{ParsevalsIdentity}
    \|\mathcal{R}(u^n)\|_{H^{-1}}^2 = \|\psi^n\|_{H_0^1}^2 \approx \sum_{k=1}^{N_{\text{test}}} |\langle \mathcal{R}(u^n), \varphi_k \rangle|^2.
\end{equation}

Using \eqref{eq:residual} then we obtain the following approximation to the dual norm of the residual and therefore a robust estimation of the error at each time step 
\begin{equation}\label{norm1time}
\begin{split}
    &\|\mathcal{R}(u^n)\|_{H^{-1}}^2 
    \approx \sum_{k=1}^{N_{\text{test}}} \Big| a(u^n, \varphi_k ) - l(\varphi_k ) \Big|^2 \\ & \quad = \sum_{k=1}^{N_{\text{test}}} \Bigg| \int_{\Omega} \mathcal{C}\, u^n\, \varphi_k \, dx  + \Delta t \int_{\Omega} \mathcal{K}\, \nabla u^n \cdot \nabla \varphi_k \, dx  - \int_{\Omega} f^n \, \varphi_k \, dx \Bigg|^2.
\end{split}
\end{equation}

Therefore, the total loss function is defined as the cumulative dual norm of the residual over all time steps. Substituting \eqref{norm1time}, we arrive at:
{\small \begin{equation}\label{LossCalor}
    \mathcal{L}(u_{\theta}) = \Delta t \sum_{n=1}^{N_{\text{time}}}  \sum_{k=1}^{N_{\text{test}}} \Bigg| \int_{\Omega} \mathcal{C}\, u^n\, \varphi_k \, dx + \Delta t \int_{\Omega} \mathcal{K}\, \nabla u^n \cdot \nabla \varphi_k \, dx  - \int_{\Omega} f^n\, \varphi_k \, dx \Bigg|^2
\end{equation}}
Minimizing \eqref{LossCalor} effectively minimizes the error of the neural network approximation across the entire space-time domain.

\begin{remark}
    For instance, if $\Omega = (a, b)$, one possible choice of orthonormal basis with respect to the inner product $(u, v)_{H_0^1} = \int_{\Omega} \nabla u \cdot \nabla v \, dx$ is
\begin{equation*}
    \varphi_k(x) = \sqrt{\frac{2}{b-a}} \cdot \sin\left[ k \left( \frac{\pi(x-a)}{b-a} \right) \right], \quad \text{with} \quad k \in \mathbb{N}.
\end{equation*}
    We remark that other choices of basis functions are possible; however, they require inverting the Gram matrix associated with the basis function representations, which is a computationally expensive operation depending on the choice of functions \cite{rojas2024robust}.
\end{remark}

\section{Numerical Experiments}
\label{sec:experiments}

In this section, we evaluate the performance of the proposed Variational Physics-Informed Neural Network (VPINN) framework through two distinct numerical experiments. 

First, we consider a linear \textit{toy problem}: the standard heat equation with constant coefficients. Since this problem possesses a known analytical solution, it serves as a benchmark to validate the accuracy of our method. 
Second, we apply the validated framework to the target industrial application: the non-linear freezing of coffee extracts. This second case involves temperature-dependent thermophysical properties, presenting a more complex challenge that demonstrates the robustness of the proposed time-stepping architecture.

\subsection{Toy Problem: Linear Heat Equation}

We first address the standard parabolic problem defined on the domain $\Omega = (0, \pi)$ and with $t_{\text{end}}=1$. We seek a weak solution to the strong formulation \eqref{eq:strong_form} with $\mathcal{C}(\cdot) = \mathcal{K}(\cdot) \equiv 1$. Here, the source term and the initial conditions are such that the analytical solution to the problem is: $$ u^*(x,t) = e^{-t}\sin(x)\cos(x/2).$$

The approximate solution is parameterized by a fully-connected feed-forward neural network implemented in \textit{TensorFlow}. The architecture consists of $5$ hidden layers with $32$ neurons per layer, utilizing the hyperbolic tangent (\(\tanh\)) as the non-linear activation function.

To enforce the homogeneous Dirichlet boundary conditions exactly, we employ the cutoff function strategy described in Section \ref{sec:DFR0}. The raw network output is multiplied by the function \( \chi(x) = x(\pi - x) \), ensuring that the trial solution vanishes at \( x=0 \) and \( x=\pi \) by construction.

The network parameters are optimized during $20^4$ iterations, using the {Adam} optimizer with a learning rate of \(10^{-2}\) with exponential decay. The variational loss function is computed using the following discretization:
\begin{itemize}

    \item {Time:} We recall that the output layer size corresponds to the number of time steps in the temporal discretization, we select \( N_{\text{time}} = 128 \) steps. 
    \item {Test Space:} We utilize the first \( N_{\text{test}} = 20 \) Fourier modes to approximate the test space \( H_0^1(\Omega) \).
    \item {Quadrature:} Spatial integrals are approximated using a midpoint rule with \( N_{\text{int}} = 128 \) randomly distributed points.
\end{itemize}

\begin{figure}[h]
    \centering
    \includegraphics[width=0.5\textwidth]{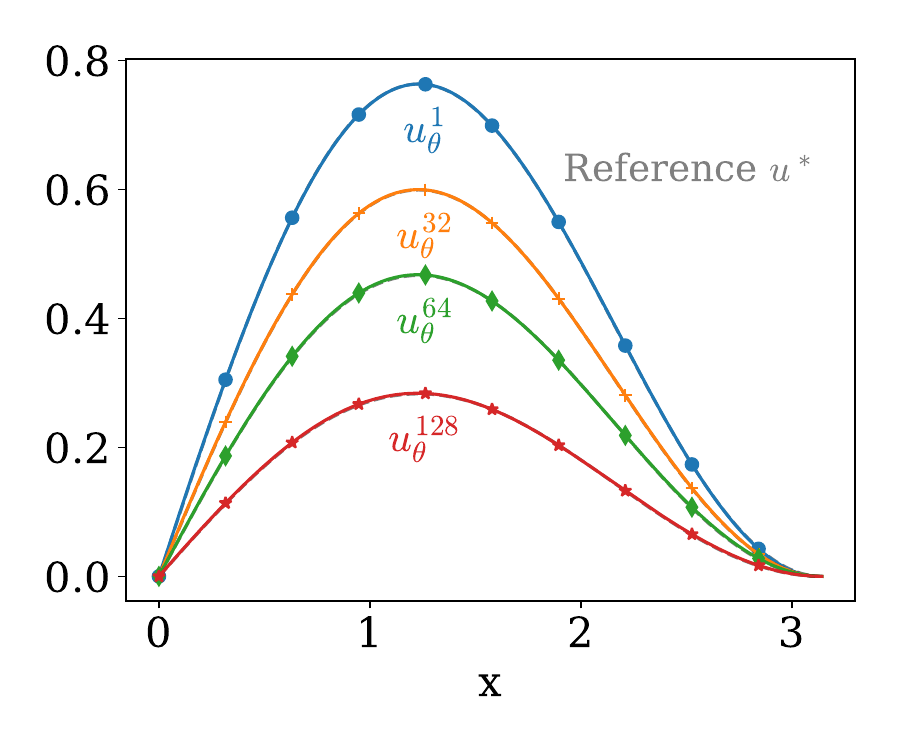}
    \caption{Comparison of the neural network approximation ($u^n_{\theta}$) and the exact solution at four different time steps (\( n=1, 32, 64, 128 \)).}
    \label{fig:solution_vs_exactPARABOLIC}
\end{figure}

Figure~\ref{fig:solution_vs_exactPARABOLIC} compares the neural network approximation \( u_{\theta} \) with the exact solution  at four distinct time steps (\( n=1, 32, 64, 128 \)). The results demonstrate that the method captures the diffusive decay of the sine wave while strictly satisfying the homogeneous boundary conditions.

\begin{figure}[htpb!]
    \centering
    \includegraphics[width=0.45\linewidth]{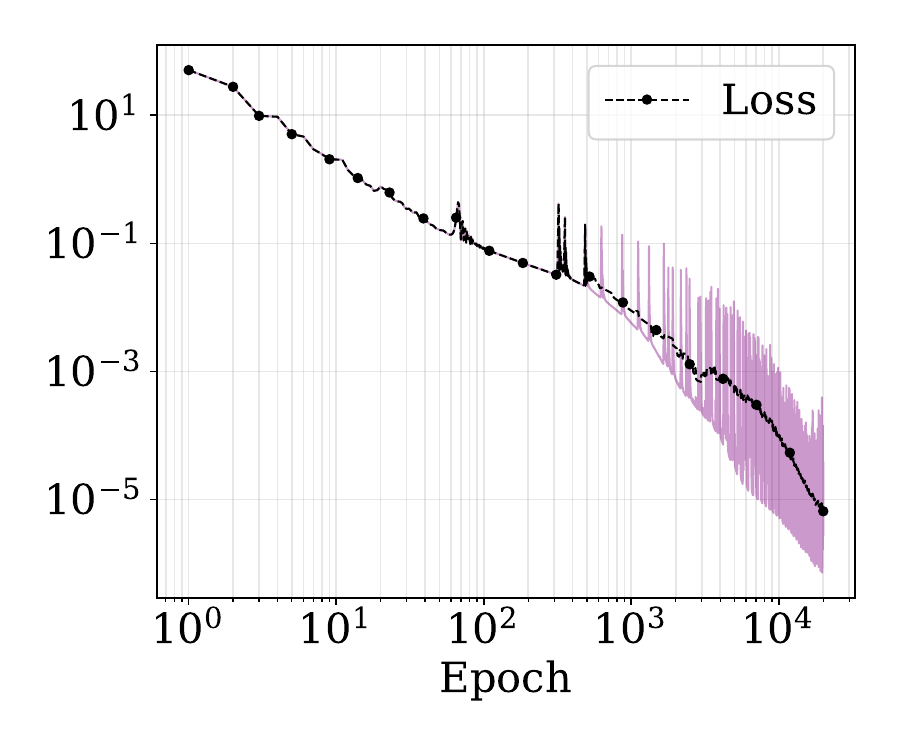} \hspace{0.2cm}
    \includegraphics[width=0.45\linewidth]{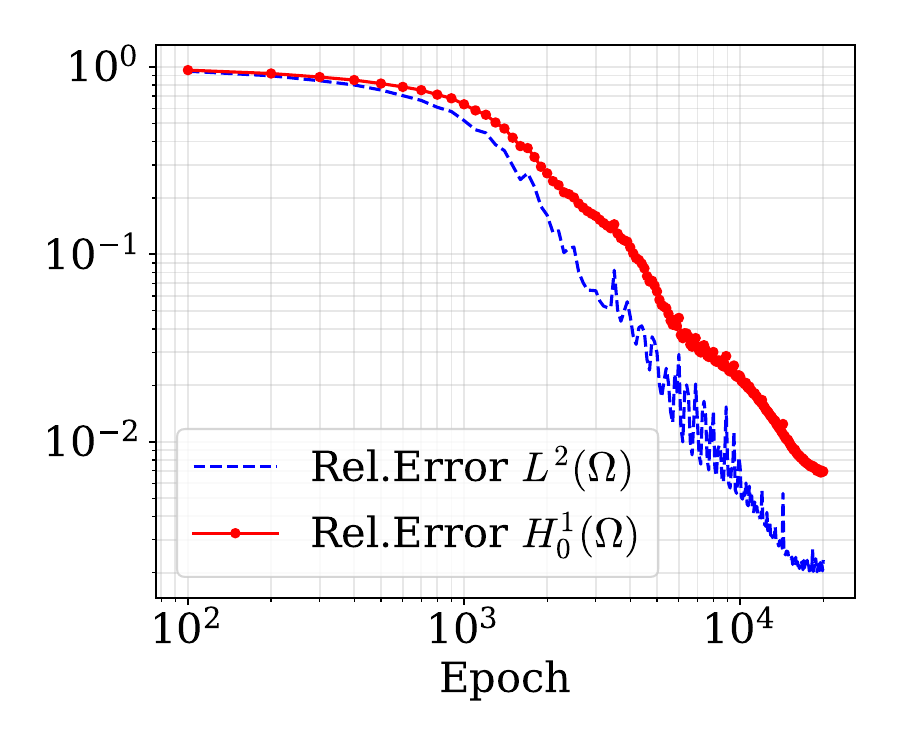}
    \caption{Evolution of the loss function during training (left). The relative errors ($L^2(\Omega)$ and $H^1_0(\Omega)$) (right).}
    \label{fig:ex1loss}
\end{figure}

Figure~\ref{fig:ex1loss} tracks the evolution of the loss function and the relative errors with respect to the exact solution $u^*(x,t)$. The reference norms are given by:
\begin{equation}
    \|u^*\|_{L^2}^2 = \int_0^T \int_{\Omega} |u^*(x,t)|^2 \, dx \, dt \quad \text{and} \quad \|u^*\|_{H^1_0}^2 = \int_0^T \int_{\Omega} \left| \frac{\partial u^*}{\partial x}(x,t) \right|^2 \, dx \, dt.
\end{equation}

These reference values ($\|u^*\|_{L^2}$ and $\|u^*\|_{H^1_0}$) are pre-computed once and used as denominators to normalize the error metrics.
In Figure~\ref{fig:ex1loss}, we observe oscillations in the loss function during convergence, attributed to the discretization limits of the variational formulation. Specifically, the test space uses a truncated basis, and the integral constraints are numerically approximated, introducing integration noise. Nevertheless, this simple example demonstrates that the numerical error decreases concurrently with the training loss, confirming that minimizing \eqref{LossCalor} effectively reduces the error (up to a constant) as expected from \eqref{eq:error_estimator}.

\subsection{Non-linear freezing of coffee extracts
}
We validate the proposed framework using data and parameters derived from the experimental setup published in \cite{restrepo2025container}. In the physical experiment, $25$ kg of coffee extract with a $31.1\%$ total dissolved solids (TDS) concentration was deposited into a cylindrical carbon-steel container. The container, equipped with an internal layered plastic bag, was subjected to a freezing process inside a forced-convection freezer maintained at $-25^\circ\text{C}$. We extract the highly temperature-dependent material properties: density, heat capacity, and thermal conductivity—directly from the physical and compositional models formulated in this experimental study.

For the numerical experiment designed to evaluate our PINN framework, we construct a 1D simplification of this system. Although the physical freezing process occurs within a cylindrical geometry, for the sake of simplicity we assume a planar 1D Cartesian domain with a characteristic domain length $d$ and assume that the domain is sufficiently large to neglect transverse effects. 

The heat transfer is governed by the one-dimensional heat conduction equation with temperature-dependent thermophysical properties:
\begin{equation}
    \begin{aligned}\label{eq:ParabolicCoffee}
       \frac{\partial}{\partial \tau} \Big( \rho(T) c_p(T) \, T \Big) - \frac{\partial}{\partial s}  \left( k(T) \frac{\partial T}{\partial s} \right) &= 0,  & &s\in (0,d) \text{ and } \tau \geq 0,\\
       T(0,\tau) = T(d,\tau) &= T_b,  &&\tau \geq 0 \\
       T(s,0) &= T_I,  & &s\in [0,d].
\end{aligned}
\end{equation}

\begin{figure}[h]
    \centering
    \includegraphics[width=\textwidth]{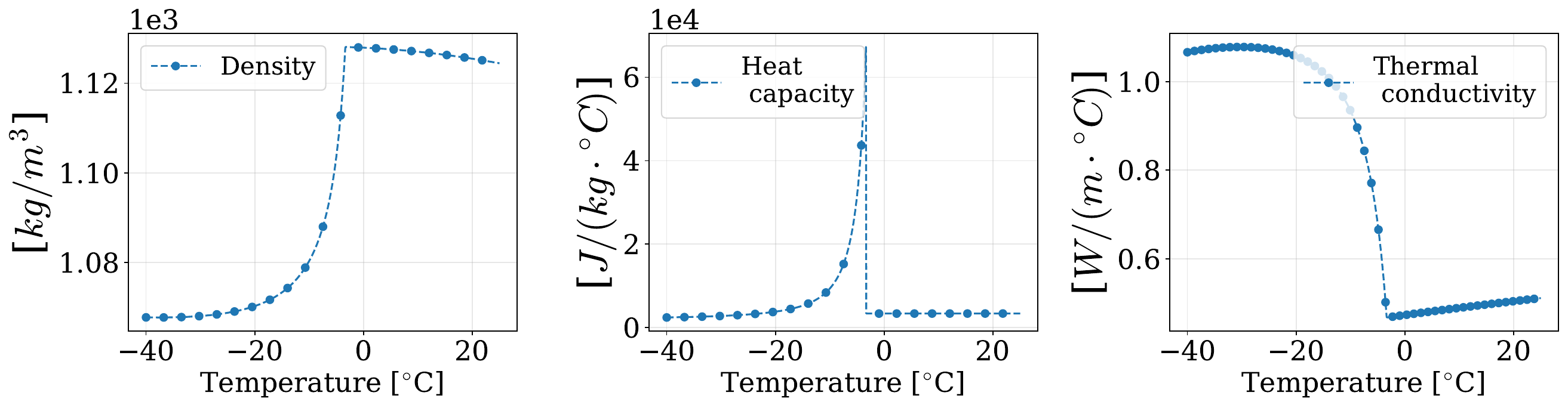}
    \caption{Thermophysical properties of the coffee extract in terms of the temperature.}
\label{fig:Thermophysical_Parameters}
\end{figure}

Here, $\tau$ denotes time $[\text{s}]$, $s$ the spatial coordinate $[\text{m}]$, and $T$ the temperature $[^\circ\text{C}]$. The parameters $\rho$, $c_p$, and $k$ represent density $[\text{kg}/\text{m}^3]$, heat capacity $[\text{J} / (\text{kg} \cdot ^\circ\text{C})]$, and thermal conductivity $[\text{W}/ (\text{m} \cdot ^\circ\text{C})]$, respectively (see Figure~\ref{fig:Thermophysical_Parameters}).

To facilitate the numerical stability and training efficiency, we introduce the following dimensionless variables:
\begin{equation}
    u = \frac{T}{T_{\text{ref}}}, \quad x = \frac{s}{d}, \quad \text{and} \quad t = \frac{\tau}{t_{\text{0}}},
\end{equation}
where $T_{\text{ref}}$ is a reference temperature (we take $T_{\text{ref}} = T_I$) and $t_{\text{0}}$ is the characteristic diffusion time scale defined as $ \displaystyle t_0 = \frac{d^2 \rho_{\text{ref}} c_{p,\text{ref}}}{k_{\text{ref}}}$, with $\rho_{\text{ref}}$, $c_{p,\text{ref}}$, and $k_{\text{ref}}$ being the reference density, specific heat, and thermal conductivity, respectively. 

Substituting these into the governing equation \eqref{eq:ParabolicCoffee} yields the dimensionless form:
\begin{equation}
     \frac{\partial}{\partial t} \Big( \mathcal{C}(u)\, u \Big) - \frac{\partial}{\partial x} \left( \mathcal{K}(u)  \frac{\partial u}{\partial x} \right) = 0, \quad x \in (0, 1), \; t > 0,
\end{equation}
where the dimensionless effective heat capacity and thermal conductivity are defined as $\mathcal{C}(u) = \frac{\rho(T) c_p(T)}{\rho_{\text{ref}} c_{p,\text{ref}}}$ and $\mathcal{K}(u) = \frac{k(T)}{k_{\text{ref}}}$, respectively.

The corresponding dimensionless boundary and initial conditions are:
\begin{equation}
\begin{aligned}
    u(0, t) = u(1, t) &= \frac{T_b}{T_I},  &&t \ge 0,  \\
    u(x, 0) &= 1,  && x \in [0, 1]. 
\end{aligned}
\end{equation}

To align with the experiments in \cite{restrepo2025container}, we set the initial temperature to $T_0 = 20^\circ\text{C}$ and impose the measured boundary temperatures of the cylinder as data-driven boundary conditions.

We parameterize the approximate solution using a fully-connected feed-forward neural network. The architecture comprises $5$ hidden layers with $32$ neurons each and utilizes the hyperbolic tangent ($\tanh$) activation function.

We optimize the network parameters over $10^5$ iterations using the Adam optimizer, applying a cosine decay learning rate with an initial value of $10^{-3}$. We compute the variational loss over a temporal discretization of $128$ time steps, which strictly matches the size of the network's output layer.

We remark that we incorporate empirical boundary data strongly in the NN using a cut-off function as explained in the previous experiment. Moreover, in this experiment we define our orthonormal basis functions in the Hilbert space $H^1(\Omega)$ rather than $H^1_0(\Omega)$. We employ $N_{\text{test}} = 64$ Fourier modes and evaluate the weak form integrals via a stochastic midpoint quadrature rule over $256$ uniformly sampled random points.

\begin{figure}[htpb!]
    \centering
    \includegraphics[width=1\textwidth]{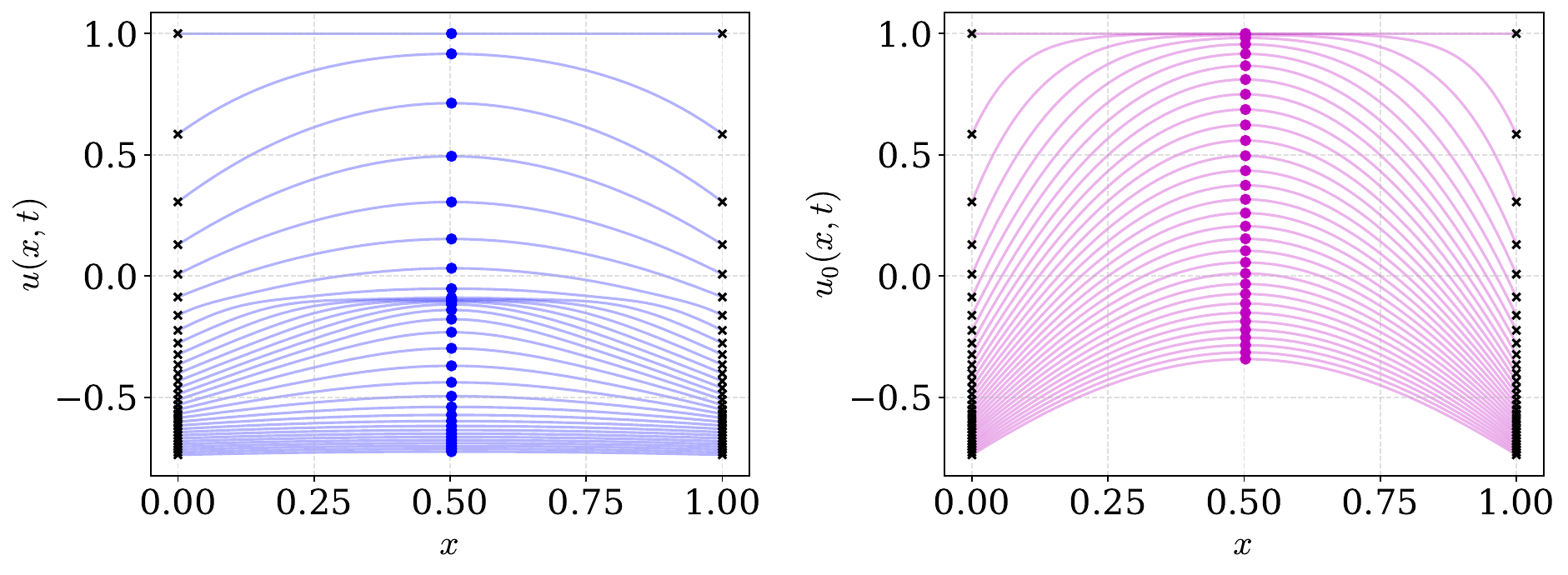} 
    \caption{Temporal evolution of the temperature profile. The non-linear solution $u$ (left) incorporates empirical thermophysical parameters, while the linear control solution $u_0$ (right) assumes constant properties. Data-driven boundary values are marked with $\times$, and the midpoint temperature is highlighted with $\bullet$.}
    \label{fig:solutionCafe}
\end{figure}

Figure~\ref{fig:solutionCafe} illustrates the temporal evolution of the temperature profiles for two scenarios:
\begin{itemize}
    \item \textbf{Non-linear model ($u$):} Incorporates the experimental, temperature-dependent thermophysical parameters presented in Figure~\ref{fig:Thermophysical_Parameters}. We emphasize that the neural network approximates the solution by evaluating the fully non-linear problem directly. This completely avoids any prior linearization of the governing equations, processing the non-linear terms natively within the loss formulation, analogous to a forward explicit scheme.
    \item \textbf{Linear model ($u_0$):} Assumes constant material parameters ($\mathcal{C} = \mathcal{K} = 1$).
\end{itemize}

\begin{figure}[htpb!]
    \centering
    \includegraphics[width=0.45\textwidth]{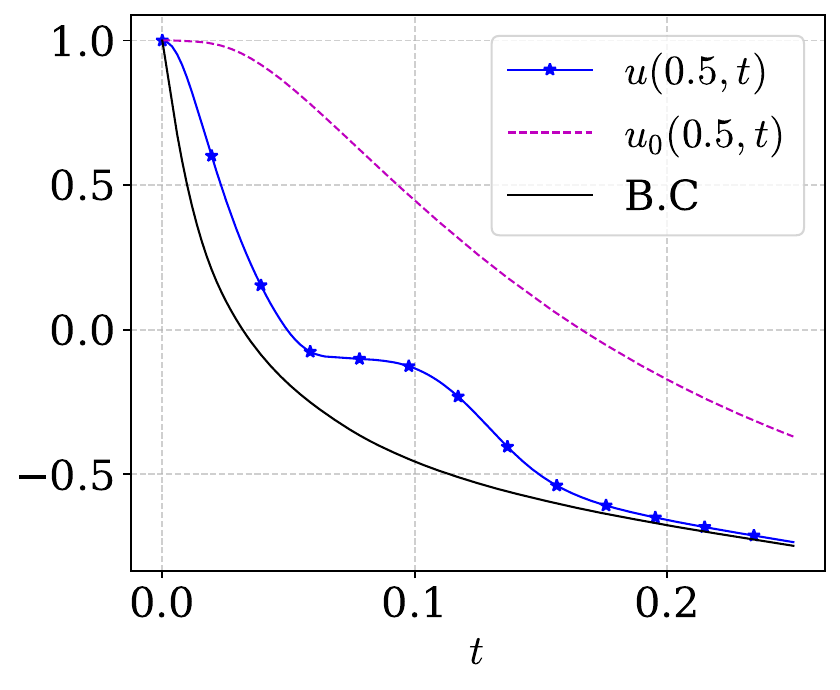}
    \caption{Evolution of the midpoint temperature for both models, plotted alongside the experimental boundary conditions.}
    \label{fig:solutionCafe2}
\end{figure}

Comparing the two solutions reveals the explicit physical impact of the non-linear coefficients. While the linear model $u_0$ diffuses at a constant rate, the non-linear solution $u$ exhibits a markedly non-uniform diffusion dynamics driven by the temperature dependence of both the effective heat capacity $\mathcal{C}(u)$ and the conductivity $\mathcal{K}(u)$. In particular, as the temperature approaches the phase-change region, the variations in thermophysical properties alter the local diffusion time scale, producing a slower thermal penetration compared to the linear case. This effect is consistent with the expected physical behavior of food freezing processes, where latent heat phenomena and structural changes in the material modify the apparent thermal inertia. The discrepancy between the linear and non-linear models is not merely quantitative but qualitative. In the linear case, the temperature profiles retain a symmetric and smoothly diffusive character throughout time. In contrast, the non-linear model shows steeper gradients near the boundaries during intermediate times, reflecting the interaction between the imposed experimental boundary data and the evolving material coefficients. This confirms that neglecting temperature-dependent properties may lead to systematic misrepresentation of transient thermal regimes in industrial applications. This difference is further highlighted in Figure~\ref{fig:solutionCafe2}, where the midpoint temperature evolution clearly separates the two scenarios. The non-linear model predicts a delayed cooling dynamics relative to the linear approximation, capturing the experimentally consistent thermal buffering effect induced by the variable heat capacity. Importantly, the VPINN framework remains stable and accurately resolves this non-linearity without any prior linearization or operator splitting. The network processes the fully non-linear weak formulation directly within the loss functional, demonstrating that the residual dual norm minimization provides sufficient structure to handle complex coefficient interactions.

\begin{figure}[htpb!]
    \centering
    \includegraphics[width=0.45\linewidth]{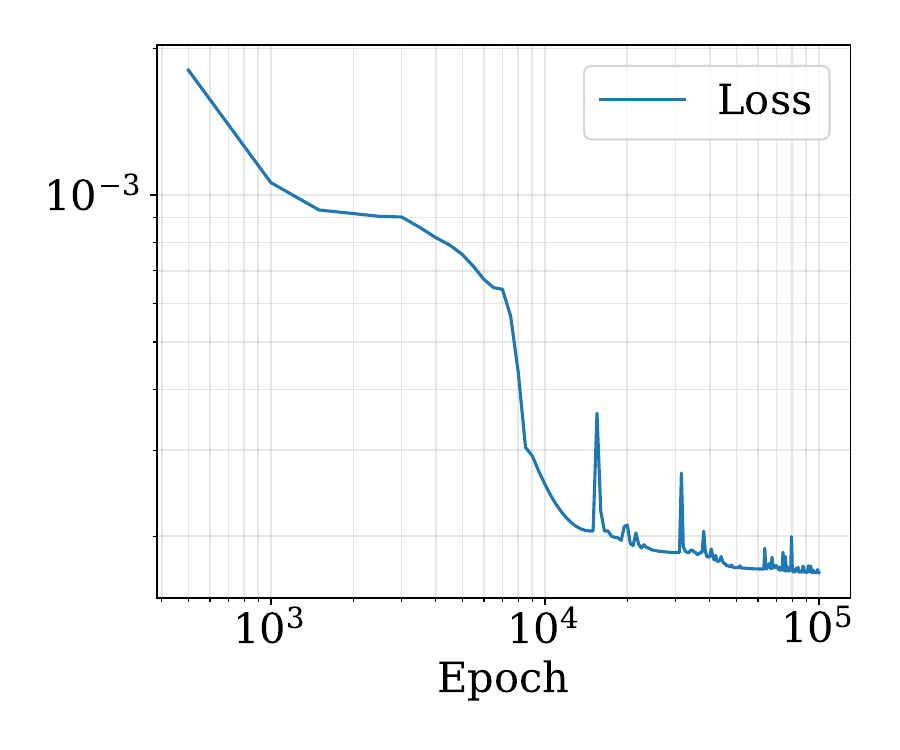}
    \caption{Evolution of the loss function over $10^5$ training iterations.}
    \label{fig:ex2loss}
\end{figure}

Finally, Figure~\ref{fig:ex2loss} displays the convergence of the loss function during training. We terminate the optimization at $10^5$ iterations, as further improvements to the neural network are not relevant. 
Overall, the results confirm that the proposed RVPINN framework is capable of solving fully non-linear parabolic problems with temperature-dependent coefficients in an industrially relevant setting, preserving physical interpretability and numerical robustness.

\section{Conclusions}
\label{sec:conclusions}
We have introduced a time-discrete Variational Physics-Informed Neural Network (RVPINN) for solving parabolic partial differential equations. The method combines classical time discretization with a residual dual norm minimization strategy, providing a theoretically grounded alternative to standard mean-squared-error PINNs. 

The numerical experiments confirm both accuracy and robustness. In particular, the industrial freezing application demonstrates that the method successfully handles temperature-dependent coefficients without linearization, capturing the physically relevant differences between linear and nonlinear diffusion dynamics.

Overall, the proposed RVPINN constitutes a mathematically consistent and computationally stable alternative to classical PINNs for transient diffusion processes. Future work will address higher-dimensional configurations and more complex multiphysics couplings.

The source code and experimental data supporting this research are available in a public repository. This includes the TensorFlow-based RVPINN implementation and all notebooks necessary to replicate the training process, figures, and analysis. The repository also hosts the empirical datasets for the measured boundary conditions and temperature-dependent thermophysical parameters. The complete repository is freely accessible at 

\url{https://github.com/manubastidas/parabolicRVPINNs}.

\section*{Acknowledgments}
This work was funded by the Universidad Nacional de Colombia, Sede Medellín, under Projects HERMES 63334 and HERMES 63067.

\section*{Declaration of Generative AI and AI-assisted technologies in the writing process}

During the preparation of this work, the author(s) used Gemini (Google) in order to assist with writing, correct English language usage, and improve text clarity and organization. After using this tool/service, the author(s) reviewed and edited the content as needed and take(s) full responsibility for the content of the publication.

\bibliographystyle{siam}
\bibliography{bibliografia}

\end{document}